\documentclass[12pt,leqno]{amsart}
\usepackage{amsmath,stmaryrd}
\usepackage{amssymb}
\usepackage{amsthm}
\usepackage{epsf}
\usepackage{graphicx}
\usepackage{esint} 
\usepackage{dsfont}
\usepackage[pagebackref]{hyperref} 
\hypersetup{pdfpagemode=FullScreen, colorlinks=true, citecolor={blue}} 
\usepackage{enumerate} 

\usepackage[authoryear]{natbib}

\numberwithin{equation}{section}

\DeclareMathOperator{\diver}{div}

\DeclareMathOperator{\supp}{supp}

\def\div{\mathop{\operatorname{div}}}

\newcommand{\ms}{\medskip}

\newcommand{\R}{\mathbb R}

\newcommand{\bp}{\noindent {\em Proof: }}
\newcommand{\ep}{\hfill $\square$ \medskip}

\newcommand{\A}{\mathcal A}

\newcommand{\B}{\mathcal B}
\newcommand{\cS}{\mathcal S}
\newcommand{\cN}{\mathcal N}

\def\YYint#1#2#3{{\setbox0=\hbox{$#1{#2#3}{\iint}$}
		\vcenter{\hbox{$#2#3$}}\kern-.51\wd0}}
 
\theoremstyle{plain}
\newtheorem{theorem}[equation]{Theorem}

\newtheorem{proposition}[equation]{Proposition}

\newtheorem{definition}[equation]{Definition}

\theoremstyle{definition}

\theoremstyle{remark}
\newtheorem{remark}[equation]{Remark}

\begin{document}

\title[The regularity problem for DKP operators]{An alternative proof of the solvability of the $L^p$-regularity problem for Dahlberg-Kenig-Pipher operators on $\R^n_+$.}

\author[Feneuil]{Joseph Feneuil}\footnote{Joseph Feneuil. Laboratoire de math\'ematiques d’Orsay, Universit\'e Paris-Saclay,
91405 Orsay, France \\
{\bf Email:} joseph.feneuil@universite-paris-saclay.fr. {\bf Orcid:} 0000-0001-5505-4450}
\address{Joseph Feneuil. Laboratoire de math\'ematiques d’Orsay, Universit\'e Paris-Saclay,
91405 Orsay, France}
\email{joseph.feneuil@universite-paris-saclay.fr}

\begin{abstract} 
In this article, we present a simpler and alternative proof of the solvability of the regularity problem - that is, the Dirichlet problem with boundary data in $\dot W^{1,p}$ - for uniformly elliptic operators on $\mathbb{R}^n_+$ under a (possibly large) Carleson measure condition. In addition, we slightly expand the class of operators for which the regularity problem is solvable, and establish an analogous result for weighted uniformly elliptic operators on $\mathbb{R}^n \setminus \mathbb{R}^d$, where $d < n - 1$.
\end{abstract}

\maketitle

\ms\noindent{\bf Keywords:} Dirichlet regularity problem, Dahlberg-Kenig-Pipher operators, Carleson measure condition, Lipschitz domains. 

\ms\noindent {\bf AMS classification:} 35J25 (primary), 31B25, 35J70.

\tableofcontents

\section{Introduction}
\label{S1}

Boundary value problems continue to attract significant attention, with one key direction of research being the identification of optimal conditions on the operator and the domain that guarantee the problem's solvability. When the boundary data lies in $L^p$, the solvability of the Dirichlet problem for the Laplacian - along with the related concept of absolute continuity of the harmonic measure - has been extensively studied; see, for example, the works of \cite{Dah77, HM14, HMU14, Azz, AHMMT}. These articles ultimately demonstrate that the Dirichlet problem is solvable in $L^p$ for some $p\in (1,\infty)$ whenever the boundary of $\Omega$ is uniformly rectifiable and $\Omega$ satisfies certain interior connectedness conditions. We will not introduce all the relevant definitions here and instead refer the reader to \cite{DS1, DS2} for the notion of uniform rectifiability, and to \cite{AHMMT} for the optimal geometric framework that guarantees the solvability of the Dirichlet problem.

Beyond the Laplacian, the Dirichlet problem has also been considered for more general uniformly elliptic operators. In this context, there are known counterexamples (see \cite{MM81, CFK81}), as well as positive results for $t$-independent operators (\cite{JK84, Ver, KKPT00, AA11, HKMP15}) and for operators satisfying a Carleson measure condition (\cite{FKP91, KP01, DPP07, DP19, DHM21}).

\smallskip

The regularity problem in $L^p$ concerns estimating the gradient of the solution in terms of the gradient of the boundary data. In fact, it can be viewed as a Dirichlet problem with boundary data in the homogeneous Sobolev space $\dot W^{1,p}$. This problem arises naturally after establishing the solvability of the classical Dirichlet problem in $L^p$, reflecting the familiar principle from the smooth setting that a higher regularity of the boundary data leads to greater regularity of the solution. For uniformly elliptic operators of the form $L:=-\diver A \nabla$ with real coefficients, the regularity problem is known to be stronger than the Dirichlet problem: solvability of the regularity problem for $L$ in $L^p$ implies solvability of the Dirichlet problem for the adjoint operator $L^*:=-\diver A^T \nabla$ in the dual space $L^{p'}$; see \cite[Theorem A.2]{MT??} for the most general version of this result. For brevity, we will denote this duality result as (R$_p$)$_L$ $\implies$ (D$_{p'}$)$_{L^*}$. Consequently, to establish the solvability of the regularity problem, it suffices to prove the converse implication (D$_{p'}$)$_{L^*}$ $\implies$ (R$_p$)$_L$.

However, a recent result shows that this converse fails in full generality: \cite[Proposition 1.9]{Gal25} provides a counterexample where (D$_{p'}$)$_{L^*}$ holds but (R$_p$)$_L$ fails when no assumptions are made about the operator or the domain. Despite this, several positive results are known. The regularity problem is solvable for $t$-independent operators in the upper half-space $\R^n_+$, as shown in \cite{JK84, KR08, HKMP15b}, and for operators satisfying a Carleson measure condition, as in \cite{KP93, KP95, DPR17, MPT??, DHP??}. Further contributions to the study of the regularity problem can be found in \cite{Shen07, DK12, DFMpert, DFMreg, GMT??, MZ??}. A particularly notable advance is due to \cite{MT??}, who established the implication (D$_{p'}$)$_{-\Delta}$ $\implies$ (R$_p$)$_{-\Delta}$ for the Laplacian in domains with uniformly rectifiable boundaries in all dimensions, and who highlighted the need to use Haj\l asz derivatives instead of classical tangential derivatives for the boundary data. This breakthrough removed a long-standing restriction that had confined most earlier proofs to Lipschitz domains.

As in the case of domains with uniformly rectifiable boundaries, the solvability of the regularity problem in $L^p$ for operators satisfying a Carleson measure condition - referred here as DKP operators\footnote{We refer to these as DKP operators, in reference to Dahlberg, who first conjectured their suitability for the Dirichlet problem, and Kenig and Pipher, who later validated the conjecture.} - remained an open problem for over two decades. Since the foundational work \cite{KP01} on the Dirichlet problem for such operators, it was conjectured that no smallness assumption on the Carleson norm was also needed for the solvability of the regularity problem. This conjecture was recently resolved simultaneously and independently in \cite{MPT??} and \cite{DHP??}. The aim of this article is fairly modest: we provide a simpler and shorter proof of the main result in \cite{DHP??}.  Our approach shows that the solvability of the regularity problem follows from a straightforward combination of two existing results: the treatment of Carleson perturbations in \cite{KP95}, and the change of variable introduced in \cite{FenDKP}. In addition, we explain why the results remain valid under a slightly weaker definition of DKP operators - one that does not require taking the supremum over Whitney balls. This refinement was already observed in the study of the Dirichlet problem by \cite{DFMdahlberg, BTZ, FLM}, and it turns out to work for the regularity problem as well. Finally, we extend the result to the higher co-dimensional setting by proving the analogous statement for operators defined on $\R^n\setminus \R^d$, which in particular generalizes the result of \cite{DFMreg}.

\subsection{Precise Statements} Let $L:=-\diver A \nabla$ be a uniformly elliptic operator on $\R^n_+$, where the coefficient matrix $A(X) \in \R^{n\times n}$ satisfies:
\begin{itemize}
	\item A boundedness condition:
    \begin{equation} \label{defboundedop}
    |A(X)\xi \cdot \zeta| \leq C |\xi||\zeta| \qquad \text{for all } X \in \mathbb{R}^n_+, \ \xi, \zeta \in \mathbb{R}^n;
    \end{equation}

    \item An ellipticity condition:
    \begin{equation} \label{defellipticop}
    A(X)\xi \cdot \xi \geq C^{-1} |\xi|^2 \qquad \text{for all } X \in \mathbb{R}^n_+, \ \xi \in \mathbb{R}^n.
    \end{equation}
\end{itemize}
In both cases, the constant $C:=C(A)>0$ does not depend on $X$. For such operators, one can define the associated elliptic measure $\{\omega_L^X\}_{X\in \R^n_+}$, a family of Borel probability measures on $\R^{n-1} = \partial \R^n_+$, characterized as follows: Given $f\in C^\infty_0(\R^{n-1})$, define
\begin{equation} \label{defuf}
u_f(X) := \int_{\mathbb{R}^{n-1}} f(y) \, d\omega_L^X(y), \qquad X \in \mathbb{R}^n_+.
\end{equation}
Then:
\begin{itemize}
    \item $u_f\in W^{1,2}(\R^n_+)$,

    \item $u_f$ is a weak solution to $Lu = 0$, i.e.,
    \begin{equation} \label{defweaksol}
    \int_{\mathbb{R}^n_+} A(X)\nabla u_f(X) \cdot \nabla v(X) \, dX = 0 \qquad \text{for all } v \in C^\infty_0(\mathbb{R}^n_+),
    \end{equation}

    \item $u_f$ extends continuously to $\R^n_+$, with the trace satisfying $u_f(x,0) = f(x)$ for all $x\in \R^{n-1}$.
\end{itemize}
Intuitively, the elliptic measure $\omega^X_L$ can be seen as the fundamental object that allows the construction of the solution of the continuous Dirichlet problem:
\[\left\{ \begin{array}{l} Lu  = 0 \text{ in } \R^n_+, \\ u = f \in C^0_0(\R^{n-1}) \text{ on } \partial \R^n_+, \\ 
u \in C^0_b(\overline{\R^n_+}).
\end{array}\right.\]

The definitions of the Dirichlet and regularity problems in $L^p$ rely on the notion of non-tangential maximal functions, which we now introduce. For each $x\in \R^{n-1}$ and a function $f$ defined in $\R^n_+$, we define the non-tangential cone
\[\Gamma(x):= \{(y,t) \in \R^n_+, \, |y-x|<t\}.\]
Using this cone, we define two functionals: the non-tangential maximal function
\[\cN(f)(x) := \sup_{X \in \Gamma(x)} |f(X)|,\]
and its averaged variant
\[\widetilde \cN(f)(x) := \sup_{(y,t) \in \Gamma(x)} \left(\fint_{B(y,2t)} \int_t^{2t} |f(z,s)|^2\, \frac{ds}{s}\, dz\right)^\frac12.\]
Alternative definitions of $\cN$ and $\widetilde \cN$ exist. For instance, the aperture of the cone $\Gamma(x)$ may be widened or narrowed, or the averaging region in $\widetilde \cN$ can be replaced by other comparable sets. These variations yield equivalent functionals in the sense that their $L^p$-norms are comparable for any $p\in (1,\infty)$, see for instance \cite[Lemma 2.1]{MPT??}, and is the only aspect relevant for our purposes.

\begin{definition} \label{DefDpRq}
We say that the Dirichlet problem for the operator $L=-\diver A \nabla$ on $\R^{n}_+$ is solvable in $L^p$ - (D$_p$)$_L$ in short - if there exists $C>0$ such that for any $f\in C^\infty_0(\R^{n-1})$, the solution $u_f$ constructed in \eqref{defuf} satisfies
\begin{equation} \label{defDp}
\|\cN(u_f)\|_{L^p(\R^{n-1})} \leq C \|f\|_{L^p(\R^{n-1})}.
\end{equation}
Moreover, we say that the regularity problem for $L$ is solvable in $L^q$ - (R$_q$)$_L$ in short -  if there exists $C>0$ such that for any $f\in C^\infty_0(\R^{n-1})$, the solution $u_f$ constructed in \eqref{defuf} satisfies
\begin{equation} \label{defRq}
\|\widetilde \cN(\nabla u_f)\|_{L^q(\R^{n-1})} \leq C \|\nabla f\|_{L^q(\R^{n-1})}.
\end{equation}
\end{definition}

Some comments on the above definitions are in order:
\begin{enumerate}
\item In \eqref{defRq}, the term $\nabla u_f$ is a full gradient in $\R^n$ while $\nabla f$ refers to the tangential  gradient in $\R^{n-1}$. This highlights a key difficulty of the regularity problem: we are attempting to control $n$ functions on $\R^n_+$ using only $n-1$ functions defined on the boundary $\R^{n-1}$. In contrast, the Dirichlet problem requires controlling just one function on $\R^n_+$ via one function on the boundary. 
\item As the name suggests, solvability of the Dirichlet problem in $L^p$, as defined in Definition \ref{DefDpRq}, implies that for any boundary data in $f \in L^p(\R^{n-1})$, one can construct a solution $u_f$ by approximating $f$ by smooth  and compactly supported functions $f_n$, and taking the limit of the corresponding $u_{f_n}$ defined via in \eqref{defuf}. While the resulting function $u_f$ may not be continuous up to the boundary if $f$ is not continuous, it still recovers $f$ in the non-tangential sense:
\[\lim_{X \in \Gamma(x)} u_f(x) = f(x) \qquad \text{ for a.e. } x\in \R^{n-1}.\]
This is classical result, see for instance \citep{FKP91}.
\item A similar density argument applies to the regularity problem. However, some care is needed because the natural space for boundary data $\{f \in L^1_{loc}(\R^{n-1}), \, \|\nabla f\|_{L^q} < \infty\}$ is an homogenous space and thus equipped only with a seminorm. For a detailed discussion of this point, we refer the reader to \citep[Section 6]{KP93}.
\end{enumerate}

Let us recall a few fundamental results concerning the Dirichlet and regularity problems, to complement the background that the reader should keep in mind while reading this article.

\begin{theorem}[{\citep[Remark 1.7.5, Theorem 1.8.14, Lemma 1.8.7]{KenigBook}}] \label{Thbasic}
Let $L=-\diver A \nabla$ be a uniformly elliptic operator on $\R^n_+$, and let $p,q\in (1,\infty)$.
\begin{enumerate}[(a)]
\item If the Dirichlet problem for $L$ is solvable in $L^{p}$, then there exists $\epsilon>0$ such that the Dirichlet problem for $L$ is solvable in $L^r$ for all $r\in (p-\epsilon,\infty)$, i.e.
\[ \text{(D$_p$)$_L$ $\implies$ (D$_r$)$_L$ \, for $r\in (p-\epsilon,\infty)$.}\]
\item If the regularity problem for $L$ is solvable in $L^{q}$, then there exists $\epsilon>0$ such that the regularity problem for $L$ is solvable in $L^q$ for all $q\in (1,q_0+\epsilon)$, i.e.
\[ \text{(R$_q$)$_L$ $\implies$ (R$_r$)$_L$ \ for $r\in (1,q+\epsilon)$.}\]
\item If the regularity problem for $L$ is solvable in $L^{q}$, then the Dirichlet problem for $L^* = -\diver A^T \nabla$ is solvable in $L^{q'}$, where $\frac1q + \frac1{q'} = 1$, i.e.
\[ \text{(R$_q$)$_L$ $\implies$ (D$_{q'}$)$_{L^*}$.}\]
\end{enumerate}
\end{theorem} 

The final point of the theorem was discussed earlier in the introduction and has been established in various geometric settings, notably in \cite[Theorem 5.4]{KP93}, \cite[Theorem 1.5]{DFMpert}, and \cite[Theorems 1.6 and A.2]{MT??}.
It is worth noting that the duality implication (R$_q$)$_L$ $\implies$ (D$_{q'}$)$_{L^*}$ fails in full generality when $L$ has complex coefficients, as shown in \cite{May10}. However, the implication can still hold for $t$-independent operators with complex coefficients, provided that suitable De Giorgi-Nash-Moser estimates are satisfied; see \cite{AM14, HKMP15b}.

The converse implication - whether the solvability of the Dirichlet problem implies the solvability of the regularity problem - is false in general. A counterexample is given in \cite[Proposition 1.9]{Gal25}. Nevertheless, under additional assumptions, the converse can be valid. For example:
\begin{itemize}
\item  If the operator is $t$-independent, the implication holds; see \cite{HKMP15b}.
\item If one knows {\em a priori} that the regularity problem is solvable for some $q_0 \in (1,\infty)$, then the full range of exponents $q$ for which the regularity problem is solvable can be determined via duality: it corresponds to the H\"older conjugates of those $p$ for which the Dirichlet problem for $L^*$ is solvable in $L^p$; see \cite{Shen07}. 
\item If the operator satisfies a Carleson measure condition and the boundary is uniformly rectifiable, the implication also holds; see \cite{MPT??}.
\end{itemize}

\medskip

The operators under consideration in our article are described in terms of Carleson measure. First, we say that $f\in CM$ - or $f\in CM(M)$ if we want to highlight the constant - if the function $f$ in $\R^{n}_+ = \{(x,t) \in \R^{n-1} \times (0,\infty)\}$ is such that $f^2 dx dt/t$ is a Carleson measure, meaning that 
\[\sup_{(z,r)\in \R^{n}_+} \fint_{B(z,r)} \int_0^r |f(x,t)|^2 \, \frac{dt}{t} \, dx \leq M^2.\]

\begin{definition}[Weak DKP operators]
We say that $L=-\diver A \nabla$ is a DKP operator if $L$ is uniformly elliptic and there exists $M>0$ such that $t|\nabla A| \in CM(M)$.

Moreover, we say that $L=-\diver A \nabla$ is a weak-DKP operator if $L$ is uniformly elliptic and if there exists $M>0$ such that
\[ \sup_{(z,r) \in \R^n_+} \fint_{B(z,r)} \int_0^r \left( \inf_{A_0 \text{ constant}} \fint_{B(x,2t)} \int_t^{2t} |A(y,s) - A_0|^2 \frac{ds}{s} dy \right) \frac{dt}{t} dx \leq M^2,\]
that is if
\[ f_L(x,t) :=  \inf_{A_0 \text{ constant}} \fint_{B(x,2t)} \int_t^{2t} |A(y,s) - A_0|^2 \frac{ds}{s} dy\]
satisfies $f_L\in CM(M)$.
\end{definition}

Weak-DKP operators can be understood as the uniformly elliptic operators whose $L^2$ oscillations over Whitney regions are controlled by a Carleson packing condition. A closely related but slightly stronger condition - used in \cite{DHP??} - requires that
\begin{equation} \label{frtg}
\inf_{A_0 \text{ constant}} \sup_{B(x,2t) \times (t,2t)} |A(y,s) - A_0| \in \text{CM}(M),
\end{equation}
mirroring the Carleson-type assumptions introduced in \cite{DPP07}. One of the main goals of this article is to show how the solvability of the regularity problem for operators satisfying \eqref{frtg} can be extended to the broader class of weak-DKP operators.

The solvability of the Dirichlet problem for weak-DKP operators is well-established; see for instance \cite[Theorem 6.9]{BTZ} and \cite[Theorem 1.21]{FLM}, as well as earlier results in \cite{KP01, DPP07, DFMdahlberg} for earlier results with more restrictive assumptions. While these results - and those in Theorem~\ref{Thbasic} - remain valid for more general domains than $\R^n_+$, we restrict our attention to the upper half-space throughout this article for simplicity and clarity of exposition.

\begin{theorem} \label{ThDpDKP}
Let $L:=-\diver A \nabla$ be a weak-DKP operator on $\R^n_+$. Then there exists $p\in (1,\infty)$ such that the Dirichlet problem for $L$ is solvable in $L^p$.
\end{theorem}

The main result of this article is the counterpart of the above theorem, now addressing the regularity problem.

\begin{theorem} \label{mainTh}
Let $L:=-\diver A \nabla$ be a weak-DKP operator on $\R^n_+$. Then there exists $p\in (1,\infty)$ such that the regularity problem for $L$ is solvable in $L^p$.
\end{theorem}

Let us emphasize that our main result is closely related to those in \citep{DHP??} and \citep{MPT??}. The goal of this article is threefold: first, to present an alternative proof of the main theorem in \citep{DHP??} that is both shorter and more transparent; second, to explain how the results of \citep{DHP??} and \citep{MPT??} - which rely on a slightly stronger condition - can be extended to the broader class of weak-DKP operators; third, we outline how our approach naturally extends to the setting of $\R^n \setminus \R^d$ (with $n < d - 1$) and to the weighted uniformly elliptic operators introduced in \cite{DFMprelim}.

\medskip

The structure of the paper is as follows. In the first part, we collect and refine several known results that are either already available in the literature or require only minor adjustments. These refinements, although modest, may be of independent interest to readers working in the area. The second part contains the proof of Theorem \ref{mainTh}, which we condense into just two pages. We then present the analogue of Theorem \ref{mainTh} in higher codimension, and explain how our methods extend seamlessly to that setting.

\medskip

In the rest of the article, we use the notation $A \lesssim B$ when there exists a constant $C>0$ independent of the relevant parameters such that $A \leq C B$, and we write $A \approx B$ when $A \lesssim B$ and $B \lesssim A$.

\medskip

\noindent {\bf Acknowledgements:} The author is grateful to the referee for their valuable comments, which have significantly enhanced the clarity and readability of the article. 

\section{Preliminaries}

Note that many of the results presented in this section - namely Theorems \ref{ThPoisson} and \ref{ThCarlpert}, and Propositions \ref{propCarleson} and \ref{propstabchg} - remain valid in significantly more general settings than the half-space $\R^n_+$. However, we choose to restrict our discussion to $\R^n_+$ in order to avoid the additional technical definitions that arise in more general geometries.

\subsection{Solvability of the Poisson-Dirichlet problem} 

We define the homogeneous space $\dot W^{1,2}(\R^n_+)$ as
\[ \dot W^{1,2}(\R^n_+) := \{ u \in L^1_{loc}(\R^n_+), \, \nabla f \in L^2(\Omega)\}.\] 
equipped with the semi-norm $\|.\|_W:=\|\nabla .\|_{L^2}$. The completion of $C^\infty_0(\R^n_+)$ under the norm $\|.\|_W$ defines a normed space denoted by $W^{1,2}_0(\R^n_+)$, which can be characterized as the subspace of $\dot W^{1,2}(\R^n_+)$ consisting of functions with vanishing trace on $\R^{n-1} = \partial \R^n_+$. The aforementioned properties of $\dot W^{1,2}(\R^n_+)$ and $W^{1,2}_0(\R^n_+)$ are classical, but for proofs and a detailed presentation of those properties (for domains far more general than $\R^n_+$), we refer the interested reader to \cite[Sections 4, 6 and 9]{DFMmixed}.

\begin{definition} \label{defPoisson}
If $L:=-\diver A \nabla$ and $\mathbf h \in L^2(\R^n_+)$, we say that $v \in \dot W^{1,2}(\R^n_+)$ solves the inhomogeneous Dirichlet problem
\[\left\{ \begin{array}{ll} L v = - \diver \mathbf h & \text{ in } \R^n_+\\ v = 0 & \text{ on } \R^{n-1} = \partial \R^n_+  \end{array}\right.\]
if $v$ satisfies the conditions:
\begin{itemize}
\item $v\in W_0^{1,2}(\R^n_+)$;
\item For any $\varphi \in W_0^{1,2}(\R^n_+)$, we have
\begin{equation} \label{defLu=divh}
 \iint_{\R^n_+} A \nabla v \cdot \nabla \varphi \, dx\, dt = \int_{\R^n_+} \mathbf{h}\cdot \nabla \varphi \, dx\, dt. 
 \end{equation}
\end{itemize}
\end{definition}

To state our next theorem, we introduce additional tools: first the truncated averaged non-tangential maximal function:
\begin{equation} \label{defS}
\cN_{p,K}(v)(x) := \sup_{(y,r) \in \Gamma(x) \atop (y,r) \in K}  \left(\fint_{|Y-(y,r)| < r/4} |v(Y)|^p \, dY\right)^\frac1p \qquad \text{ for } K \Subset \R^n_+,
\end{equation}
and then the square functions $\A$ and $\cS$:
\begin{equation} \label{defA}
\A(v)(x) := \left(\int_{\Gamma(x)} v^2 \frac{dt\, dx}{t^{n}}\right)^\frac12 \quad \text{ and } \quad \cS(v)(x) := \left(\int_{\Gamma(x)} |t\nabla v|^2 \frac{dt\, dx}{t^{n}}\right)^\frac12,
\end{equation}
meaning that $\cS(v) = \mathcal A(t\nabla v)$.

\begin{theorem} \label{ThPoisson}
Let $q\in (1,\infty)$. For any compact set $K\Subset \R^n_+$ and any $\mathbf{F} \in L^1_{loc}(\R^n_+,\R^n)$, there exists a compactly supported function $\mathbf h \in L^\infty(\R^n_+)$ satisfying
\begin{equation} \label{eqduality}
\|\cN_{1,K}(\mathbf{F})\|_{L^q(\R^{n-1})} \leq 2 \iint_{\R^n_+} \mathbf{F}\cdot \mathbf{h}\, dt\, dx.
\end{equation}

Moreover, we can choose $\mathbf h$ so that the following holds: if $L := - \div A \nabla$ is a uniformly elliptic operator on $\R^{n}_+$ for which the Dirichlet problem for $L^*$ is solvable in $L^{q'}$, then there exists a $C>0$ independent of $\mathbf{F}$ and $K$ (and hence $\mathbf h$) such that the solution  $v \in \dot W^{1,2}(\R^n_+)$ to the inhomogeneous Dirichlet problem
\[\left\{ \begin{array}{ll} L^* v = - \diver \mathbf h & \text{ in } \R^n_+\\ v = 0 & \text{ on } \R^{n-1} = \partial \R^n_+  \end{array}\right.\]
satisfies the estimate
\begin{equation} \label{S<N<h}
\|\cS(v)\|_{L^{q'}(\R^{n-1})}  + \|\cN(v)\|_{L^{q'}(\R^{n-1})} \leq C.
\end{equation}
\end{theorem}

\begin{remark}
The estimate on $\|\cN(v)\|_{L^{q'}(\R^{n-1})}$ in \eqref{S<N<h} is closely related to the Poisson-Dirichlet problem studied in \cite{MPT??}. However, the framework developed in \cite{MPT??} does not provide bounds on the square function $\cS(v)$, which is the quantity that matters to us. For this reason, we will instead follow the original approach of Kenig and Pipher in \cite{KP95},  which was subsequently extended to more general settings in \cite{DFMpert}, and use it to establish the desired estimate.
\end{remark}

\bp The proof of the theorem can be traced back to \cite{KP95}, which treats the case where the domain is a ball and the operator is self-adjoint, and to \cite{DFMpert}, which extends the argument to unbounded domains and operators that are not necessarily self-adjoint.  However, in both references, the result appears as a combination of intermediate steps. For clarity and completeness, we explain here how our statement follows from \cite{DFMpert}.

The choice of $\mathbf h$ is made in \cite[Lemma 4.1]{DFMpert}. There, a functional $T: L^\infty_{loc}(\R^n_+) \to L^1_{loc}(\R^{n-1})$ is introduced (see equation (4.13) in \cite{DFMpert}) as an auxiliary tool. According to \cite[Lemma 4.4]{DFMpert}, this functional is such that, for the given $\mathbf h$, we have
\begin{equation} \label{pol2}
 \|T(\mathbf h)\|_{L^{q'}} \leq  C
 \end{equation}
where $C$ depends only on $q'$ and $n$.

Since $\mathbf h$ is bounded and compactly supported (in particular,  $\mathbf h \in L^2$), the function $v$ can be represented using the Green function $G_{L}$ associated with the operator $L$ as follows
\[v(X):= \int_{\R^n_+} \nabla_Y G_L(Y,X) \cdot \mathbf h(Y)\, dY \quad \text{ for a.e. } X\in \R^{n}_+.\]
Using standard estimates on the Green function, together with \cite[Corollary 4.8 and Lemma 4.10]{DFMpert}, we obtain, for all $p\in (0,\infty)$,
\begin{equation} \label{pol1}
\|\cS(v)\|_{L^p(\R^{n-1})} + \|\cN(v)\|_{L^p(\R^{n-1})}  \leq C_p \|\mathcal M_{\omega_{L^*}}(T(\mathbf h))\|_{L^p(\R^{n-1})},
\end{equation}
where $C_p>$ depend only on $p$, the dimension $n$, and the ellipticity and boundedness constants of $L$. Here $\omega_{L^*}$ is the elliptic measure with pole at infinity associated to $L^*$, and $\mathcal M_{\omega_{L^*}}$ is the Hardy-Littlewood maximal function with respect to that measure. 
The solvability of the Dirichlet problem for $L^*$ in $L^{q'}$ is equivalent to the boundedness of $\mathcal M_{\omega_{L^*}}$ on $L^{q'}(\R^{n-1},dx)$, see for instance Theorem 4.13 in \cite{KenigBook}. Therefore, the bound \eqref{pol1} becomes
\[\|\cS(v)\|_{L^{q'}(\R^{n-1})} + \|\cN(v)\|_{L^{q'}(\R^{n-1})}  \lesssim \|T(\mathbf h)\|_{L^{q'}(\R^{n-1})} \lesssim 1\]
by \eqref{pol2}. The theorem follows.
\ep

\subsection{Carleson inequality}

\begin{proposition} \label{propCarleson}
There exists $C>0$ such that for any $a\in CM(M)$ and any couple $(f,g)$ of functions on $\R^n_+$, we have
\[ \iint_{\R^n_+} a(x,t) f(x,t) g(x,t) \, dx\, \frac{dt}{t} \leq CM \int_{\R^{n-1}} \cN(f)(x) \A(g)(x) \, dx.\]
\end{proposition}

\bp
The proof relies on a classical stopping time argument commonly used in such estimates - see, for instance, the proof of (3.33) in \citep{DHP??}. Alternatively, one can adapt the original argument from Coifman, Meyer and Stein (\citep[Theorem 1]{CMS85}) to our situation with few modifications.
\ep

\subsection{Carleson perturbations}

\begin{theorem} \label{ThCarlpert}
Let $L_0:= - \diver A_0 \nabla $ and $L_1 := -\diver A_1 \nabla$ be two uniformly elliptic operators defined on $\R^n_+$. Assume that there exists $C>0$ such that 
\begin{equation} \label{tNA<C}
\sup_{(x,t) \in \R^n_+} |t\nabla A_0(x,t)| \leq C
\end{equation}
and $|A_0 - A_1| \in CM(M)$, that is
\begin{equation} \label{A0-A1carl}
\sup_{x,r \in \R^n_+} \fint_{|x-y|<r} \int_0^r |A_0(y,t) - A_1(y,t)|^2 \frac{dt}{t} \, dx \leq M^2.
\end{equation}
Assume moreover that the regularity problem for $L_0$ is solvable in $L^p$ for some $p\in (1,\infty)$. Then there exists $q \in (1,p)$ such that the regularity problem is solvable in $L^q$.
\end{theorem}

\bp
The differences between our approach and those found in the existing literature (see \citep{KP95}) are as follows:
\begin{itemize}
\item We impose the additional assumption \eqref{tNA<C}.
\item Our result assumes the $L^2$-Carleson perturbation condition \eqref{A0-A1carl}, which is strictly weaker than the standard $L^\infty$-Carleson perturbation condition:
\begin{equation} \label{A0-A1carl2}
\sup_{x,r \in \R^n_+} \fint_{|x-y|<r} \int_0^r \left(\sup_{|(z,s)-(x,t)|<t/2} |A_0(z,s) - A_1(z,s)|^2\right) \frac{dt}{t} , dx \leq M^2.
\end{equation}
\end{itemize}

Let $f\in C^\infty_0(\R^{n-1})$, and let $u_{0,f}$ and $u_{1,f}$ denote the solutions of the Dirichlet problem with the boundary data $f$ for the operators $L_0$ and $L_1$, respectively. Let $\mathbf F := u_0 - u_1$ and let $\mathbf h$ defined from $\mathbf F$ as in Theorem \ref{ThPoisson}. Since the regularity problem is solvable in $L^p$ for $L_0$, it follows from Theorem \ref{Thbasic} that the Dirichlet problem for $L_0^*$ is solvable in $L^{p'}$. Because $L_1^*$ is a Carleson perturbation of $L_0^*$, there exists $q\in (1,p]$ such that the Dirichlet problem for $L_1$ is solvable in $L^{q'}$  (see \cite[Corollary 1.16]{Fen??} for $L^2$-Carleson perturbations of the Dirichlet problem, and \citep[Theorem 2.3]{FKP91}--\citep[Corollary 1.31]{FP22} for earlier results on $L^\infty$-Carleson perturbations). Since $q<p$, we have the solvability of the regularity problem for $L_0$ in $L^q$ (see \citep[Theorem 5.2]{KP93}), that is
\[\|\widetilde \cN (\nabla u_{0,f})\|_q \leq C \|\nabla f\|_q.\]
Moreover using the assumption \eqref{tNA<C}, we have the following Moser-type estimate
\begin{equation} \label{MoserNu}
|\nabla u_{0,f}(x,r)| \leq C \left( \fint_{B(x,r)} \int_{r/3}^{3r} |\nabla u_{0,f}(y,t)|^2 \, \frac{dt}{t}\, dy \right)^\frac12.
\end{equation}
This follows from Lemma 3.1 in \citep{GW82}, combined with the Moser estimate and the Poincar\'e inequality: indeed, we have
\begin{multline*}
|\nabla u_{0,f}(x,r)| \leq \frac{C}{r} \sup_{(y,t) \in B(x,r/2) \times (r/2,2r)} |u_{0,f}(y,t) - u_{0,f,x,r}| \\ \lesssim  \frac{1}{r} \left(\fint_{B(x,r)} \int_{r/3}^{3r} |u_{0,f}(y,t) - u_{0,f,x,r}|^2 \, \frac{dt}{t}\, dy\right)^\frac12 \\ \lesssim \left(\fint_{B(x,r)} \int_{r/3}^{3r} |\nabla u_{0,f}(y,t)|^2 \, \frac{dt}{t}\, dy\right)^\frac12,
\end{multline*}
where $u_{0,f,x,r}$ is the average of $u_{0,f}$ on $B(x,r) \times (r/3,3r)$. From there, we deduce $\cN(\nabla u_{0,f}) \lesssim \widetilde \cN(\nabla u_{0,f})$, hence
\begin{equation} \label{NNu0<Nf}
\|\cN(\nabla u_{0,f})\|_q \leq C \|\nabla f\|_q.
\end{equation}

By the definition of $\mathbf h$, we have
\begin{equation} \label{identityFh=Nu0Nv}
\|\cN_{1,K}(\nabla F)\|_q \leq 2 \iint_{\R^n_+} \nabla F \cdot \mathbf h \, dxdt = 2 \iint_{\R^n_+} (A_0-A_1) \nabla u_{0,f} \cdot \nabla v\, dxdt,
\end{equation}
where $v\in W^{1,2}_{0}(\R^n_+)$ is the solution to $-\diver A_1^* \nabla v = -\diver \mathbf h$ as in Theorem \ref{ThPoisson}. The identity \eqref{identityFh=Nu0Nv} follows as in \citep[Corollary 2.9]{KP95} or \citep[Lemma 3.18]{CHM19}. Applying the Carleson inequality (Proposition \ref{propCarleson}) together with the H\"older inequality yields:
\[\|\cN_{1,K}(\nabla F)\|_q \lesssim M \|\cN(\nabla u_{0,f})\|_q \|S(v)\|_{L^{q'}}
\lesssim M \|\cN(\nabla u_{0,f})\|_q\]
using estimate \eqref{S<N<h} for the last bound. Since $u_{1,f}$ is a solution to $L_1u_{1,f} = 0$, the gradient $\nabla u_{1,f}$ satisfies a reverse $L^2$-$L^1$ H\"older inequality. Consequently,
\[\|\cN_{2,K}(\nabla u_{1,f})\|_q \leq C \|\cN_{1,K}(\nabla u_{1,f})\|_q \lesssim (1+M) \|\cN(\nabla u_{0,f})\|_q \lesssim  (1+M)  \|\nabla f\|_q\]
by \eqref{NNu0<Nf}. The constant being independent of $K \Subset \R^n_+$, we take $K \uparrow \R^n_+$ to conclude.
\ep

\subsection{Change of variable}

\begin{proposition} \label{propwDKP}
Assume that $L=-\diver A \nabla$ is a weak DKP operator. Then there exists a constant $M>0$, a function $M' : \epsilon \in (0,1) \to (1,\infty)$ such that, for any $\epsilon>0$, we have a decomposition $A(x,t) = B(x,t) + C(x,t)$ satisfying
\begin{enumerate}[(i)]
\item the matrix $B$ is elliptic with the same elliptic constant as $A$,
\item $\sup_{(x,t) \in \R^n_+} |t\nabla B(x,t)| \leq \epsilon$,
\item $t \nabla B \in CM(M)$,
\item $C \in CM(M'(\epsilon))$.
\end{enumerate}
\end{proposition}


\bp
Assume first that the Proposition is true for some large $\epsilon_0$, which depends only on $n$ and the ellipticity constant of $L$. That is $A = B_1 + C_1$ where $B_1$ has the same (or smaller) elliptic constant as $A$, $|t\nabla B_1(x,t)| \leq \epsilon_0$, and where $B_1$ and $C_1$ satisfy $|t\nabla B_1| + |C_1| \in CM(M_0)$.

The main idea is that $B$ is a smooth average of $B_1$ and so $C=B-B_1 + C_1$. Take $\phi \in C^\infty_0(\R^{n-1})$ be such that $\phi \geq 0$, $\supp \phi \subset B(0,1)$, and $\int \phi\, dx = 1$, and $\psi \in C^\infty(\R)$ be such that $\psi\geq 0$, $\supp \psi \in (1,2)$, and $\int_\R \psi(\ln(s)) ds/s = 1$. For some $\Lambda >1$ to be chosen later, we define 
\[B_\Lambda(x,t) :=  \int_{\R^{n-1}} \int_0^\infty \underbrace{\frac{ s^{1-n}}{\ln(\Lambda)} \phi\Big( \frac{y-x}{s} \Big) \psi\Big( \frac{\ln(s/t)}{\ln(\Lambda)}\Big)}_{=:\Phi_{x,t,\Lambda}(y,s)} B_1(y,s) \, \frac{ds}{s}   \, dy.\]
Note that $B_\Lambda(x,t)$ is an average of $B_1(x,t)$, since
\begin{equation} \label{BLambdaaverage}
\int_{\R^{n-1}} \int_0^\infty \Phi_{x,t,\Lambda}(y,s) \, \frac{ds}{s}   \, dy = 1,
\end{equation}
so $(i)$ is satisfied by construction, and moreover 
\begin{equation} \label{suppPhi}
\supp \Phi_{x,t,\Lambda}(y,s) \subset W_{x,t,\Lambda}:= \big\{(y,s)\in \R^{n}_+, \, s  \in (\Lambda t,\Lambda^2t) \text{ and } y \in \B(x,s)\big\}.
\end{equation}
For some constants $c_\phi$ and $c_\psi$ that depend only on $\phi$ and $\psi$, we have
\begin{multline} \label{tNxBL}
|t\nabla_x B_\Lambda(x,t)| \leq\|B_1\|_{L^\infty(\R^n_+)} \int_{\R^{n-1}} \int_0^\infty \frac{t s^{-n}}{\ln(\Lambda)} \Big|\nabla \phi\Big( \frac{y-x}{s} \Big)\Big| \psi\Big( \frac{\ln(s/t)}{\ln(\Lambda)}\Big)  \, \frac{ds}{s}   \, dy \\
\leq \frac{\|B_1\|_{L^\infty(\R^n_+)}}{\Lambda} \int_{\R^{n-1}} \int_0^\infty \frac{s^{1-n}}{\ln(\Lambda)} \Big|\nabla \phi\Big( \frac{y-x}{s} \Big)\Big| \psi\Big( \frac{\ln(s/t)}{\ln(\Lambda)}\Big)  \, \frac{ds}{s}   \, dy = \frac{c_\phi \|B_1\|_{L^\infty(\R^n_+)}}{\Lambda}
\end{multline}
and
\begin{multline} \label{tdtBL}
|t\partial_t B_\Lambda(x,t)| \leq \frac{\|B_1\|_{L^\infty(\R^n_+)}}{\ln(\Lambda)} \int_{\R^{n-1}} \int_0^\infty \frac{s^{1-n}}{\ln(\Lambda)} \phi\Big( \frac{y-x}{s} \Big) \Big|\psi'\Big( \frac{\ln(s/t)}{\ln(\Lambda)}\Big)\Big|  \, \frac{ds}{s}   \, dy \\
= \frac{c_\psi \|B_1\|_{L^\infty(\R^n_+)}}{\ln(\Lambda)}.
\end{multline}
We take $\Lambda$ large enough (depending on $\epsilon)$ to satisfy conclusion $(ii)$. The quantity $B_\Lambda$ is morally the convolution of $B_1$ with a mollifier, so we find that $t\nabla B_\Lambda$ is a (weighted) average of $t\nabla B_1$ over the Whitney area $B(x,\Lambda^2t) \times (\Lambda t,\Lambda^2t)$. In fact, we have
\[\nabla_x B_\Lambda(x,t) = \int_{\R^{n-1}} \int_0^\infty \Phi_{x,t,\Lambda}(y,s)  \nabla_x B_1(y,s) \, \frac{ds}{s} \, dy,\]
and
\[ t\partial_t B_\Lambda(x,t) = \int_{\R^{n-1}} \int_0^\infty \Phi_{x,t,\Lambda}(y,s) \Big( s \partial_s B_1(y,s) - (y-x)\cdot \nabla_x B_1(y,s) \Big)  \, \frac{ds}{s} \, dy.\]
Thus, we obtain the bound
\begin{multline} \label{tNBL}
|t\nabla B_\Lambda(x,t)| \leq \frac{C_{\phi,\psi}}{\ln\Lambda}  \int_{\Lambda t}^{\Lambda^2 t} \fint_{B(x,s)} |s\nabla B_1(y,s)| \, dy \, \frac{ds}{s} \\
\leq C_{\phi,\psi}   \left(\frac{1}{\ln\Lambda} \int_{\Lambda t}^{\Lambda^2 t} \fint_{B(x,s)} |s\nabla B_1(y,s)|^2 \, dy \, \frac{ds}{s} \right)^\frac12
\end{multline}
by using \eqref{BLambdaaverage}, \eqref{suppPhi}, and the Jensen inequality.

Now, we check that $(iii)$ is satisfied. Specifically, we need to show that $t\nabla B_\Lambda$ satisfies the Carleson measure condition with a constant independent on $\Lambda$. We begin by writing
\begin{multline*}
\fint_{B(z,r)} \int_0^r |t\nabla B_\Lambda(x,t)|^2 \frac{dt}{t} dx = \fint_{B(z,r)} \int_0^{r/\Lambda^2} |t\nabla B_\Lambda(x,t)|^2 \frac{dt}{t} dx \\
+ \fint_{B(z,r)} \int_{r/\Lambda^2}^r |t\nabla B_\Lambda(x,t)|^2 \frac{dt}{t} dx =: I_1 + I_2.
\end{multline*}
We treat $I_2$ using \eqref{tNxBL}--\eqref{tdtBL}, and we easily get that
\[I_2 \leq C_{\phi,\psi} \left( \frac{\|B_1\|_{L^\infty(\R^n_+}}{\ln \Lambda} \right)^2 \fint_{B(z,r)} \int_{r/\Lambda^2}^r  \frac{dt}{t} dx \leq \frac{C}{\ln\Lambda}\]
for a constant $C$ that depends on $\phi$, $\psi$, and the boundedness constant of $B_1$. Thus, as $\Lambda$ increases, $I_2$ becomes smaller.  To estimate $I_1$, we use \eqref{tNBL} and Fubini's lemma to obtain 
\begin{multline*}
I_1 \leq \frac{C_{\phi,\psi}}{\ln\Lambda} \fint_{B(z,2r)} \int_0^r  |s\nabla B_1(y,s)|^2 \left( \int_{s/\Lambda^2}^s \fint_{B(y,s)} dx \, \frac{dt}{t} \right) \frac{ds}{s} \, dy \\
\leq  C_{\phi,\psi} \fint_{B(z,2r)} \int_0^r  |s\nabla B_1(y,s)|^2 \frac{ds}{s} \, dy \leq C_{\phi,\psi} M_0,
\end{multline*}
where $M_0$ is such that $t\nabla B_1(x,t) \in CM(M_0)$. Conclusion $(iii)$ follows.

Finally, proving $(iv)$ is straightforward, as $B_\Lambda$ is some average of $B_1$, and thus the Poincar\'e inequality gives
\[\int_{t}^{2t} \fint_{B(x,t)} |B_\Lambda(y,s) - B_1(y,s)|^2 \, dy \, \frac{ds}{s} \leq C_\Lambda  \int_{t}^{2\Lambda^2 t} \fint_{B(x,2\Lambda^2t)} |s \nabla B_1(y,s)|^2 \, dy \, \frac{ds}{s}\]
For a given tent set $B(z,r)\times (0,r)$, we cover it with a  finitely overlapping collection of Whitney regions $\{B(x_i,t_i)\times (t_i,2t_i)\}_{i\in I}$. Since the collection  $\{B(x_i,2\Lambda^2t_i)\times (t_i,2\Lambda^2 t_i)\}_{i\in I}$ is also finitely overlapping, we deduce that 
\[ \fint_{B(z,r)} \int_0^r |B_\Lambda(x,t) - B_1|^2 \frac{dt}{t} dx \leq C_\Lambda \fint_{B(z,4\Lambda^2 r)} \int_0^{4\Lambda^2 r} |t\nabla B_1(x,t)|^2 \frac{dt}{t} dx \leq C'_\Lambda M_0,\]
as desired.

\medskip

It remains to demonstrate that we can construct the decomposition $A = B_1 + C_1$ in the first place. The computations follow closely those in the proof of \cite[Corollary 2.3]{DPP07}. For instance, we can define $B_1(x,t)$ as
\[ B_1(x,t):= \int_{R^{n-1}} \int_0^\infty \Phi_{x,t,\Lambda}(y,s) A(y,s) \, \frac{ds}{s} \, dy\]
for $\Lambda = 2^{1/4}$. The proof of the fact that $t\nabla B_1(x,t) \in L^\infty(\R^n_+)$ is identical to the proof of \eqref{tNxBL}--\eqref{tdtBL}. To show that $|t\nabla B_1| + |C_1| \in CM(M_0)$, we proceed as follows: for any constant matrix $A_0$, observe that 
\[|t\nabla B_1(x,t)| \lesssim \fint_{B(x,2t)} \int_{t}^{2t} |A(y,s) - A_0| \, \frac{ds}{s} \, dy\]
and 
\[\fint_{B(x,2^{1/4}t)} \int_{t}^{2^{1/4}t} |C_1(y,s)| \, \frac{ds}{s} \, dy \lesssim \fint_{B(x,2t)} \int_{t}^{2t} |A(y,s) - A_0| \, \frac{ds}{s} \, dy\]
thanks to \eqref{BLambdaaverage}. The Carleson bounds on $|t\nabla B_1|$ and $C_1$ follow then from the fact that $A$ is weak-DKP. The details are left to the reader.
\ep

\begin{proposition} \label{propchgvar}
Let $L=-\diver A \nabla$ be a weak DKP operator on $\R^n_+$. Then there exists a bi-Lipschitz change of variable $\rho$ on $\R^n_+$ such that 
\begin{enumerate}[(a)]
\item $\rho(x,0) = (x,0)$ for any $x\in \R^{n-1}$,
\item The conjugate of $L$ by $\rho$ - that is the elliptic operator $L_\rho:= -\diver A_\rho \nabla$ for which $L_\rho (u\circ \rho) = 0$ whenever $Lu=0$ - can be decomposed as 
\[A_\rho = \begin{bmatrix} B_{\rho,||} & \mathbf{b}_{\rho} \\ 0 & 1 \end{bmatrix} + C_\rho,\]
where $|t\nabla B_{\rho,||}| + |t\nabla \mathbf b_{\rho}| + |C_\rho| \in CM$.
\end{enumerate} 
\end{proposition}

\bp
The goal of \citep{FenDKP} is to present this change of variables in detail, and we refer the reader to this article for the full construction.

Since $L$ is a weak DKP operator, Proposition \ref{propwDKP} ensures that we can decompose $A$ into $B + C$, where $B$ has the same elliptic constant as $A$, $\sup_{\R^n_+} |t\nabla B| \leq \epsilon$, and $|t\nabla B| + |C| \in CM(M_\epsilon)$. The parameter $\epsilon$, which may be chosen arbitrary small, will be fixed later. We write $B$ in block form as
\[B(x,t) =: \begin{bmatrix} B_{||}(x,t) & \mathbf b(x,t) \\ \mathbf v(x,t) & h(x,t)\end{bmatrix}\]
and we construct the Lipschitz map
\[\rho(x,t) := (x+ t\mathbf v(x,t),th(x,t))\]
which obviously fixes $\R^{n-1}$. For $\epsilon$ small enough depending only on the elliptic constant of $B$ (hence $A$), the function $\rho$ is invertible, which makes it a bi-Lipschitz change of variable from $\R^n_+$ onto itself. A straightforward computation - provided in \citep{FenDKP} - shows that the conjugate of $L$ under $\rho$ satisfies the statement $(b)$ of the Proposition.
\ep

Let us finish this paragraph with the fact that the solvability of the regularity problem is stable under bi-Lipschitz changes of variable.

\begin{proposition}[Proposition 3.8 in \citep{FenDKP}] \label{propstabchg}
Let $L=-\diver A \nabla$ be a uniformly elliptic operator on $\R^n_+$ and $\rho$ be a bi-Lipschitz map from $\R^n_+$ to $\R^n_+$ that fixes the boundary $\R^{n-1}$. For any $p\in (1,\infty)$, the solvability of the regularity problem for $L$ in $L^p$ is equivalent to the solvability of the regularity problem for the conjugate operator $L_\rho$ from $L$ by $\rho$. 
Moreover, the constants appearing in \eqref{defRq} for $L$ and $L_\rho$ are comparable, with a ratio depending only on the dimension $n$ and the bi-Lipschitz constant of $\rho$.
\end{proposition}

\subsection{The case of the Laplacian}

\begin{proposition} \label{propLaplacian}
The regularity problem for the Laplacian in $\R^n_+$ is solvable in $L^p$ for all $p\in (1,\infty)$, that is for any $f\in C^\infty_0(\R^{n-1})$, the harmonic extension $u_f$ of $f$ satisfies 
\[ \|\cS(\nabla u_f)\|_{L^p(\R^{n-1})} \approx \|\cN(\nabla u_f)\|_{L^p(\R^{n-1})} \leq C \|\nabla f\|_{L^p(\R^{n-1})},\]
where $u_f$ is given by convolution with the Poisson kernel in the upper half-space.
\end{proposition}

\bp
This result is classical, though the precise attribution is unclear, as it can be deduced from the expression of $u_f$ using the Poisson kernel. 

As a special case of \cite{DPP07},  we can deduce that the Dirichlet problem for the Laplacian in $\R^n_+$ is solvable in $L^p$ for all $p\in (1,\infty)$. Together with \cite{KKPT00}, which proves the equivalence between the $L^p$-norm of $\cS(u_f)$ and $\cN(u_f)$, we deduce
\[ \|\cS(u_f)\|_{L^p(\R^{n-1})} \approx \|\cN(u_f)\|_{L^p(\R^{n-1})} \leq C \|f\|_{L^p(\R^{n-1})}.\]
From there, the solvability of the regularity problem is simple. Indeed, any of the tangential derivatives of $u$ is also harmonic, so we have
\[ \|\cS(\nabla_x u_f)\|_{L^p(\R^{n-1})} \approx \|\cN(\nabla_x u_f)\|_{L^p(\R^{n-1})} \leq C \|\nabla f\|_{L^p(\R^{n-1})}.\]
As for the normal derivative of $u$, it is also harmonic, so we still have 
\[ \|\cN(\partial_t u_f)\|_{L^p(\R^{n-1})} \approx \|\cS(\partial_t u_f)\|_{L^p(\R^{n-1})} \leq \|\cS(\nabla_x u_f)\|_{L^p(\R^{n-1})} + \|\mathcal A(t\partial_{tt}^2 u_f)\|_{L^p(\R^{n-1})}.\] 
But since $u$ is harmonic, we have $|\partial_{tt}^2 u| \leq |\nabla_x^2 u|$, meaning that we even get 
\[ \|\cN(\partial_t u_f)\|_{L^p(\R^{n-1})} \approx \|\cS(\partial_t u_f)\|_{L^p(\R^{n-1})} \lesssim \|\cS(\nabla_x u_f)\|_{L^p(\R^{n-1})} \lesssim \|\nabla f\|_{L^p(\R^{n-1})}.\]
The proposition follows.
\ep

\section{Proof of the main theorem}

\subsection{Case of the codimension 1}
We now turn to the proof of Theorem \ref{mainTh}, using the results developed in the previous section. Thanks to Propositions \ref{propchgvar} and \ref{propstabchg}, it suffices to establish the $L^p$ solvability of operators $L:=-\diver A \nabla$ under the following assumptions on the coeffcient matrix $A$:
\begin{itemize}
\item $A$ can be decomposed as
\[A = \begin{bmatrix} B_{||} & \mathbf{b} \\ 0 & 1 \end{bmatrix} + C;\]
\item The decomposition verifies $t|\nabla B_{||}| + t |\nabla \mathbf b| + |C| \in CM$.
\end{itemize}
By Theorem \ref{ThCarlpert}, we can disregard the Carleson perturbation $C$, reducing the problem to studying operators of the form $L:=-\diver B \nabla$ where 
\[ B = \begin{bmatrix} B_{||} & \mathbf b \\ \mathbf{0} & 1 \end{bmatrix}.\]
Moreover, Theorem \ref{ThDpDKP} ensures that for DKP-type operators, the Dirichlet problem is solvable in $L^{p'}$ for some $p'\in (1,\infty)$.  Hence, for the class of operators described above, it remains to prove that (D$_{p'}$)$_{L^*}$ $\implies$ (R$_p$)$_L$. To summarize, we need to show the following theorem.

\begin{theorem} \label{mainTh2}
Let $L=-\diver B \nabla$ is a uniformly elliptic operator on $\R^n_+$ satisfying $|t\nabla B| \in CM$ and 
\[ B = \begin{bmatrix} B_{||} & \mathbf b \\ \mathbf{0} & 1 \end{bmatrix}.\]
Let $p\in (1,\infty)$ be such that Dirichlet problem for adjoint operator $L^* = - \diver B^T \nabla$  is solvable in $L^{p'}$, where $\frac1p + \frac1{p'}=1$. Then the regularity problem for $L$ is solvable in $L^p$. 
\end{theorem}

\bp
The proof follows a strategy similar to that of Theorem \ref{ThCarlpert}. Given any $f\in C^\infty_0(\R^{n-1})$ and any compact set $K \Subset \R^n_+$, we aim to establish the localized estimate
\begin{equation} \label{claimMainTh}
\|\cN_{1,K}(\nabla u_f)\|_{L^p} \leq C\|\nabla f\|_{L^p},
\end{equation}
where $u_f$ is the solution to $Lu_f =0$ with trace $f$ defined via the elliptic measure, and where $C$ is a constant independent of $f$ and $K$. Then by using a $L^2$-$L^1$ reverse H\"older estimate on the gradient of elliptic solutions, and then passing to the limit $K \uparrow \R^{n}_+$, we conclude that 
\[\|\widetilde \cN(\nabla u_f)\|_{L^p} \leq C\|\nabla f\|_{L^p}\]
as desired.

Let $\tilde u_f$ be the harmonic extension of $f$ (given by the harmonic measure, for instance). The functions $u_f$ and $\tilde u_f$ belong to $\dot W^{1,2}(\R^n_+)$, are (H\"older) continuous on $\overline{\R^n_+}$, and satisfy $u_f(x,0) = \tilde u_f(x,0) = f(x)$. As such, $u_f - \tilde u_f \in W^{1,2}_0(\R^n_+)$.

By Theorem \ref{ThPoisson}, there exists a bounded and compactly supported vector function $\mathbf h$ such that 
\[\|\cN_{1,K}(\nabla[u_f - \tilde u_f])\|_{L^p} \leq 2 \iint_{\R^n_+}\nabla[u_f - \tilde u_f] \cdot \mathbf h \, dx\, dt.\]
Let now $v$ be the solution to $L^*v = - \div \mathbf h$ with vanishing trace, as contructed in Theorem \ref{ThPoisson}. Since the Dirichlet problem for $L^*$ is solvable in $L^{p'}$, Theorem \ref{ThPoisson} gives the estimate
\begin{equation} \label{estSv}
\|\mathcal S(v)\|_{L^{p'}} \lesssim 1.
\end{equation}
Using the weak formulation of $L^*v = -\div \mathbf h$ - i.e. equation \eqref{defLu=divh} -  and the fact that $u_f - \tilde u_f \in W^{1,2}_0(\R^n_+)$, we have
\[\iint_{\R^n_+}\nabla[u_f - \tilde u_f] \cdot \mathbf h \, dx\, dt = 2 \iint_{\R^n_+}B\nabla[u_f - \tilde u_f] \cdot \nabla v \, dx\, dt,\]
so that
\[\|\cN_{1,K}(\nabla[u_f - \tilde u_f])\|_{L^p} \leq 2 \iint_{\R^n_+} B \nabla (u_f - \tilde u_f)  \cdot  \nabla v \,  dx\, dt.\]
Since $v \in W^{1,2}_0(\R^n_+)$, we may use it as a test function in the weak formulations of both $u_f$ and $\tilde u_f$, which are solutions to $Lu_f = 0$ and $-\Delta \tilde u_f = 0$ respectively. Subtracting the two identities yields 
\[\|\cN_{1,K}(\nabla[u_f - \tilde u_f])\|_{L^p} \leq 2 \iint_{\R^n_+} (B-I) \nabla \tilde u_f \cdot \nabla v \, dx\, dt\] 

Let us denote $D:=B-I$. This matrix satisfies the structural property
\begin{equation} \label{prD1}
D_{nj} = 0 \ \text{ for $1 \leq j \leq n$} 
\end{equation}
and, for some $M>0$, we have the Carleson condition
\begin{equation} \label{prD2}
|t\nabla D_{ij}| \in CM(M) \ \text{ for $1\leq i \leq n-1$ and $1 \leq j \leq n$}.
\end{equation}
We aim to bound the quantity
\[ T_{ij} = \iint_{\R^n_+} D_{ij} (\partial_{j} \tilde u_f) (\partial_i v) \, dx \, dt \]
for $i\in \{1,\dots,n-1\}$  (corresponding to the tangential derivatives of $v$) and $j\in \{1,\dots,n\}$. For this, observe that $\partial_n t= 1$, so we may rewrite the above as
\[T_{ij} = \iint_{\R^n_+} D_{ij} (\partial_{j} \tilde u_f) (\partial_i v) (\partial_n t)\, dx \, dt.\]
We want to perform an integration by parts in the $t$-direction, aiming to shift the $\partial_n $ derivative away from $t$. We are willing to allow two derivatives on $\tilde u_f$, since we have the control over $\cS(\nabla \tilde u_f)$ from Proposition \ref{propLaplacian}. However, we can only allow at most one derivative on $v$. Therefore, when the derivative $\partial_n$ lands on $\partial_i v$, we perform a second integration by parts - this time in the tangential direction - to move $\partial_i$ away from $v$ (noting that $\partial_i t = 0$).  Following this strategy, we arrive at the identity:
\begin{multline*}
T_{ij} = - \iint_{\R^n_+} (\partial_n D_{ij}) (\partial_{j} \tilde u_f) (\partial_i v) \, t\, dx \, dt - \iint_{\R^n_+} D_{ij} (\partial_n \partial_{j} \tilde u_f) (\partial_i v) \, t\, dx \, dt \\
+ \iint_{\R^n_+} (\partial_i D_{ij}) (\partial_{j} \tilde u_f) (\partial_n v) \, t\, dx \, dt + \iint_{\R^n_+} D_{ij} (\partial_i \partial_{j} \tilde u_f) (\partial_n v) \, t\, dx \, dt \\
=: T_{ij}^1 + T_{ij}^2 + T_{ij}^3 + T_{ij}^4.
\end{multline*}
We first note that no boundary terms arise in the integration by parts. This is justified as follows: the coefficients $D_{ij}$ are bounded since $B$ (and before $A$) is bounded, and $\partial_j \tilde u_f$ is also bounded since $\tilde u_f$ is the harmonic extension of a smooth and compactly supported data. Moreover, the term $t\partial_i v$ vanishes on average over Whitney regions of $\R^n_+$ as $t\to 0$. This follows from the fact that $\mathbf h$ is compactly supported, which implies that $v$ is a H\"older continuous, vanishes at the boundary and at infinity, and thus satisfies the decay $t\nabla v \to 0$ on average in Whitney region as a consequence of the Cacciopoli inequality.

We now estimate the four resulting terms coming from the previous integration by parts. The terms $T_{ij}^1$ and $T_{ij}^3$ are bounded using the Carleson condition \eqref{prD2} on $t\nabla D_{ij}$. Indeed, we have
\begin{multline*}
|T_{ij}^1| + |T_{ij}^3| \lesssim \iint_{\R^n_+} |t \nabla  D_{ij}| | \nabla \tilde u_f| |\nabla v| \, dx \, dt \lesssim M \int_{\R^{n-1}} \cN(\nabla \tilde u_f)(x) \mathcal A(t\nabla v)(x) \, dx \\ 
\lesssim M \|\cN(\nabla \tilde u_f)\|_{L^p} \|\cS(v)\|_{L^{p'}} \lesssim M\|\nabla f\|_{L^p}
\end{multline*}
where we used Proposition \ref{propCarleson} and \eqref{prD2} in the second step, followed by the H\"older inequality, Proposition \ref{propLaplacian}, and the bound \eqref{estSv} in the final line. For the two remaining terms, we proceed similarly:
\[|T_{ij}^2| + |T_{ij}^4| \lesssim \iint_{\R^n_+} |D_{ij}| |t\nabla^2 \tilde u_f| |\nabla v| \, dx \, dt \lesssim \|S(\nabla \tilde u_f)\|_{L^p} \|\cS(v)\|_{L^{p'}} \lesssim \|\nabla f\|_{L^p}\]
where we again used H\"older's inequality, Proposition \ref{propLaplacian}, and \eqref{estSv}.

Combining all the estimates, we obtain
\[\|\cN_{1,K}(\nabla[u_f - \tilde u_f])\|_{L^p}  \lesssim (1+M) \|\nabla f\|_{L^p}\]
and therefore
\[\|\cN_{1,K}(\nabla u_f)\|_{L^p}  \leq \|\cN(\nabla \tilde u_f)\|_{L^p} + C (1+M) \|\nabla f\|_{L^p} \lesssim (1+M) \|\nabla f\|_{L^p}\]
where the final inequality uses Proposition \ref{propLaplacian} once more.
\ep

\subsection{Case $\Omega = \R^{n}\setminus \R^d$.} In this paragraph, we discuss the analogue of Theorem \ref{mainTh} for the domain $\R^n \setminus \R^d$ defined as $\{(x,t) \in \R^d \times (\R^{n-d}\setminus \{0\})\}$, where $d<n-1$. In this context, an elliptic theory has been developed in \citep{DFMprelim}. Here, uniformly elliptic operators are those that can be expressed in the form $L = - \diver [|t|^{d+1-n} A\nabla]$, where the matrix function $A$ satisfies \eqref{defboundedop}--\eqref{defellipticop}. The cones in this setting are $\Gamma(x):= \{(y,t)\in \R^n \setminus \R^d, \, |y-x|<|t|\}$. The averaged non-tangential maximal function is given by
\[\widetilde\cN(f)(x):= \sup_{(y,t) \in \Gamma(x)} \left( \fint_{B(y,2t)} \int_{|t|\leq |s| \leq |2t|} |f(z,s)|^2 \frac{ds}{|s|^{n-d}} dz \right)^\frac12.\]
The definition of $\cN$ and the definition of the solvability of the Dirichlet and the regularity problem are now straightforward.

\begin{theorem} \label{ThRn-Rd}
Let $L=-\diver[|t|^{d+1-n} A \nabla]$ be a uniformly elliptic operator such that $A$ can be decomposed as
\[ A = B + C\]
with $|C| \in CM$, and 
\begin{enumerate}[(i)]
\item either $B$ can be written\footnote{$Id_{n-d}$ is the identity matrix of order $n-d$, $b_4$ is a scalar function, which makes $B_1$, $B_2$, and $B_3$ matrices of size $d\times d$, $d\times (n-d)$, and $(n-d)\times d$ respectively} as
\[B = \begin{bmatrix} B_1 & B_2 \\ B_3 & b_4 Id_{n-d} \end{bmatrix}\]
with $|t||\nabla B| \in CM$,
\item or $B$ can be written\footnote{$t$ and $\mathbf b_3$ are horizontal vectors, and $\mathbf b_2$ is a vertical vector} as
\[B = \begin{bmatrix} B_1 & \mathbf b_2 \frac{t}{|t|} \\ \frac{t^T}{|t|} \mathbf b_3 & b_4 Id_{n-d} \end{bmatrix}\]
with $|t|(|\nabla B_1|+ |\nabla \mathbf b_2| + |\nabla \mathbf b_3| + |\nabla b_4|) \in CM$
\end{enumerate} 
Then there exists $p\in (1,\infty)$ such that the regularity problem for $L$ is solvable in $L^p$.
\end{theorem}

\begin{remark}
Case $(ii)$ involves the types of operators obtained by rotating a weak-DKP operator on $\R^{d+1}_+$ around the boundary $\R^d$. For details on this construction, see Section 4.1 in \cite{DFMKenig}. This case, particularly when the Carleson constant is small, is also examined in \cite{DFMreg}.
Case $(i)$, on the other hand, is novel, even for operators defined with small Carleson constants. It is particularly suited to bi-Lipschitz changes of variables. Specifically, after applying a change of variable to flatten a Lipschitz graph, the off-diagonal terms in $B$ do not conform to the forms $B_2 = \mathbf b_2 \frac{t}{|t|}$ and $B_3 \frac{t^T}{|t|} \mathbf b_3$ as seen in $(ii)$. However, the off-diagonal terms could satisfy $|t| (|\nabla B_2| + |\nabla B_3|) \in CM$, if, for instance, one use the change of variable from \cite{KP01}.
\end{remark}

\noindent {\em Explanations of the proof.} 

\noindent {\bf 1-Adaptation of the preliminary results.} Theorem \ref{ThPoisson} can be extended to degenerate operators and higher codimensional boundaries, since the arguments for Theorem \ref{ThPoisson} are derived from \cite{DFMpert}, which addresses these cases. The Carleson inequality (Proposition \ref{propCarleson}) is a real variable argument that can be easily adapted to our context.
The analogue of Theorem \ref{ThCarlpert} in $\R^n \setminus \R^d$ can be proven using \cite{DFMpert}, similar to how Theorem \ref{ThCarlpert} is proven using \cite{KP95}.
Proposition \ref{propwDKP} can be adapted to $\R^n \setminus \R^d$ for both scenarios $(i)$ and $(ii)$. Specifically, case $(i)$ has already been addressed in \cite{FenDKP}, while case $(ii)$ represents a simpler variant of case $(i)$ as long as Proposition \ref{propwDKP} is concerned.
Furthermore, the analogue of Proposition \ref{propLaplacian} in $\R^n \setminus \R^d$ pertains to the solvability of the regularity problem for the operator $L_0 := -\diver[|t|^{d+1-n}\nabla]$. The operator $L_0$ is essentially the``Laplacian rotated around $\R^d$", meaning that the solution to $L_0 u= 0$ with boundary data $f$ is $u_{0,f}(x,t) := u_f(x,|t|)$, where $u_f$ is the harmonic extension of $f$ to $\R^{d+1}_+$. Therefore, the solvability of the regularity problem for $L_0$ is a direct consequence of Proposition \ref{propLaplacian}.

All these preliminary results enable us to reduce Theorem \ref{ThRn-Rd} to the case where $L=-\diver[|t|^{d+1-n} B\nabla]$, with $B$ given by
\[B = \begin{bmatrix} B_1 & B_2 \\ \mathbf 0 & I_{n-d} \end{bmatrix}\]
and either $|t||\nabla B_2| \in CM$ - in Case (i) - or $B_2 = \mathbf b_2 t/|t|$ with $|t||\nabla \mathbf b_2| \in CM$ - in Case $(ii)$. 

\smallskip

\noindent {\bf 2-Adaptation of Theorem \ref{mainTh2}.} To treat this simplified case, we follow the proof of Theorem \ref{mainTh2} and we will ultimately need to estimate
\[T_{ij} = \iint_{\R^n \setminus \R^d} D_{ij} (\partial_j \tilde u_f)(\partial_i v) \, dx \, \frac{dt}{|t|^{n-d-1}} \]
where $i\in \{1,\dots,d\}$, $j\in \{1,\dots,n\}$, $D_{ij} = B_{ij}-1$, and $\tilde u_f$ is the solution to $L_0 u = 0$ with boundary data $f$.

We express the derivatives in cylindrical coordinates, which involves using the following derivatives:
\[\partial_r := \sum_{d <k \leq n} \frac{t_k}{|t|} \partial_k,\]
corresponding to the radial derivative, and for $d < \alpha,\beta \leq n$,
\[\partial_{\varphi_{\alpha\beta}} := - \frac{t_\alpha}{|t|} \partial_\beta + \frac{t_\beta}{|t|} \partial_\alpha,\]
which corresponds to the angular derivatives. Here we identify $t \in \R^{n-d}$ with the vector $(0,\dots,0,t_{d+1},\dots,t_n)\in \R^n$. Note that any $\partial_j$, $j>d$, can be written as a linear combination of $\partial_r$ and $\partial_{\varphi_{\alpha\beta}}$. Specifically, 
\[\partial_j = \frac{t_j}{|t|} \partial_r + \sum_{k>d} \frac{t_k}{|t|} \partial_{\varphi_{jk}}.\]
After verifying that $t_k/|t|$ behaves like an angle - meaning that $\partial_r(t_k/|t|) = 0$ -  and that $\tilde u_f$ is radial - i.e. $\partial_{\varphi_{\alpha\beta}} \tilde u_f = 0$ - bounding the $T_{ij}$ in both cases ($|t||\nabla B_2| \in CM$ or $B_2 = \mathbf b_2 t/|t|$ with $|t||\nabla \mathbf b_2| \in CM$) reduces to bounding terms in the form
\[U_{ij} := \iint_{\R^n\setminus \R^d} \tilde D_{ij} (\partial_j \tilde u_f) (\partial_i v) \, dx \, \frac{dt}{|t|^{n-d-1}} \qquad i,j \in \{1,\dots, d\},\]
and
\[U_{ir} := \iint_{\R^n\setminus \R^d} \widetilde D_{ir} (\partial_r \tilde u_f) (\partial_i v) \, dx \, \frac{dt}{|t|^{n-d-1}}  \qquad i\in \{1,\dots,d\}.\]
Here, $|t||\nabla_{x,r} \tilde D| \in CM$, where $\nabla_{x,r}$ denotes the gradient involving the tangential derivative $\partial_i$, $1 \leq i \leq d$, and the radial derivative $\partial_r$ (excluding the angular derivatives $\partial_{\varphi_{\alpha\beta}}$). $\tilde D$ represents either $\tilde D_{ij}$ or $\tilde D_{ir}$. The proof is then the same as before. We use the fact that $\partial_r |t| = 1$, we integrate by parts to shift $\partial_r$ away from $|t|$ to the other terms; note that the integration by parts with respect to $\partial_r$ is 
\[\iint_{\R^n\setminus \R^n} f (\partial_r |t|) \, \frac{dt}{|t|^{n-d-1}}\, dx = - \iint_{\R^n\setminus \R^n} (\partial_r f) |t| \, \frac{dt}{|t|^{n-d-1}}\, dx,\]
provided that $|t|f(t) \to 0$ as $t\to 0$ or $\infty$. We then move $\partial_i$ away from $v$ when $\partial_r$ hits $\partial_i v$. The details are the same as what we did for $T_{ij}$ in the proof of Theorem \ref{mainTh2} and are left to the reader. 
\ep

\bibliography{refs} 
\bibliographystyle{alpha}


\end{document}